\documentclass[12pt,dvips]{amsart}
\usepackage{amsfonts, amssymb, latexsym, epsfig, url}

\newtheorem{Theorem}{Theorem}[section]

\newtheorem{Definition-Proposition}[Theorem]{Definition-Theorem}
\newtheorem{Main Conjecture}[Theorem]{Main Conjecture}

\newtheorem{Conjecture}[Theorem]{Conjecture}

\theoremstyle{remark}
\newtheorem{Example}[Theorem]{Example}

\newtheorem{Definition}[Theorem]{Definition}

\newcommand\Schub{{\mathfrak S}}

\newcommand{\Groth}{\mathfrak{G}}
\newcommand{\Demazure}{\kappa}
\newcommand{\KDemazure}{\Omega}

\theoremstyle{plain}

\newcommand{\cellsize}{11}
\newlength{\cellsz} 
\newsavebox{\cell}
\sbox{\cell}{\begin{picture}(\cellsize,\cellsize)
\put(0,0){\line(1,0){\cellsize}}
\put(0,0){\line(0,1){\cellsize}}
\put(\cellsize,0){\line(0,1){\cellsize}}
\put(0,\cellsize){\line(1,0){\cellsize}}
\end{picture}}
\newcommand\cellify[1]{\def\thearg{#1}\def\nothing{}%
\ifx\thearg\nothing
\vrule width0pt height\cellsz depth0pt\else
\hbox to 0pt{\usebox{\cell} \hss}\fi%
\vbox to \cellsz{
\vss
\hbox to \cellsz{\hss$#1$\hss}
\vss}}
\newcommand\tableau[1]{\vtop{\let\\\cr
\baselineskip -16000pt \lineskiplimit 16000pt \lineskip 0pt
\ialign{&\cellify{##}\cr#1\crcr}}}
\newcommand{\kellsize}{20}
\newlength{\kellsz} 
\newsavebox{\kell}
\sbox{\kell}{\begin{picture}(\kellsize,\kellsize)
\put(0,0){\line(1,0){\kellsize}}
\put(0,0){\line(0,1){\kellsize}}
\put(\kellsize,0){\line(0,1){\kellsize}}
\put(0,\kellsize){\line(1,0){\kellsize}}
\end{picture}}
\newcommand\kellify[1]{\def\thearg{#1}\def\nothing{}%
\ifx\thearg\nothing
\vrule width0pt height\kellsz depth0pt\else
\hbox to 0pt{\usebox{\kell} \hss}\fi%
\vbox to \kellsz{
\vss
\hbox to \kellsz{\hss$#1$\hss}
\vss}}
\newcommand\ktableau[1]{\vtop{\let\\\cr
\baselineskip -16000pt \lineskiplimit 16000pt \lineskip 0pt
\ialign{&\kellify{##}\cr#1\crcr}}}

\newcommand{\sellsize}{36}
\newlength{\sellsz} 
\newsavebox{\sell}
\sbox{\sell}{\begin{picture}(\sellsize,20)
\put(0,0){\line(1,0){\sellsize}}
\put(0,0){\line(0,1){\sellsize}}
\put(\sellsize,0){\line(0,1){\sellsize}}
\put(0,\sellsize){\line(1,0){\sellsize}}
\end{picture}}
\newcommand\sellify[1]{\def\thearg{#1}\def\nothing{}%
\ifx\thearg\nothing
\vrule width0pt height\sellsz depth0pt\else
\hbox to 0pt{\usebox{\sell} \hss}\fi%
\vbox to \sellsz{
\vss
\hbox to \sellsz{\hss$#1$\hss}
\vss}}
\newcommand\stableau[1]{\vtop{\let\\\cr
\baselineskip -16000pt \lineskiplimit 16000pt \lineskip 0pt
\ialign{&\sellify{##}\cr#1\crcr}}}

\begin{document}
\pagestyle{plain}

\title{The ``Grothendieck to Lascoux'' conjecture}
\author{Victor Reiner}
\address{Dept.~of Mathematics, Univ.~of Minnesota, Minneapolis, MN 55455, USA}
\email{reiner@umn.edu}
\author{Alexander Yong}
\address{Dept.~of Mathematics, Univ.~of Illinois at
Urbana-Champaign, Urbana, IL 61801, USA}
\email{ayong@illinois.edu}

\date{February 18, 2021}

\maketitle

\begin{abstract}
This report 
formulates a conjectural combinatorial rule that positively expands
Grothendieck polynomials into Lascoux polynomials.  It generalizes
one such formula expanding Schubert polynomials into key polynomials, 
and refines another one expanding stable Grothendieck polynomials.
\end{abstract}

\section{The open problem}

We set up the notation needed to state the problem, Conjecture~\ref{conj:main} below.
\subsection{Grothendieck and Lascoux polynomials} Define operators $\partial_i, \pi_i$ 
on polynomials $f\in {\mathbb Z}[\beta][x_1,\ldots,x_n]$
\[\partial_i(f) =\frac{1-s_i f}{x_i-x_{i+1}} \mbox{\ and  \ } \pi_i(f)=\partial_i((1+\beta x_{i+1})f)\]
where $s_{i}=(i\leftrightarrow i+1)$ is a simple transposition in the symmetric group $S_n$. The transposition $s_i$ acts on 
$f\in {\mathbb Z}[\beta][x_1,\ldots,x_n]$ by
permuting $x_i$ and $x_{i+1}$. 

\begin{Definition}[{A.~Lascoux-M.-P.~Sch\"utzenberger~\cite{LasSch2}}]
The $\beta$-\emph{Grothendieck polynomial} ${\mathfrak G}_w^{(\beta)}$ is recursively defined by the initial condition
\[ {\mathfrak G}^{(\beta)}_{w_0}=x_1^{n-1}x_2^{n-2}\cdots x_{n-1}\]
where $w_0$ is the \emph{longest permutation} that swaps  $i \leftrightarrow n+1-i$, and then setting 
\[{\mathfrak G}^{(\beta)}_w=\pi_i({\mathfrak G}^{(\beta)}_{ws_i}) \text{\ if $w(i)<w(i+1)$.}\]
\end{Definition}
\noindent
The $\beta$ parameter was introduced by S.~Fomin--A.~N.~Kirillov~\cite{Fomin.Kirillov};
\cite{LasSch2} uses $\beta=-1$.

Define a further family of operators $\widetilde{\pi}_i$ on polynomials \emph{via}
\[{\widetilde \pi}_i (f)=\partial_i(x_i(1+\beta x_{i+1})f).\]
Also, let ${\sf Comp}$ be the set of \emph{(weak) compositions}, that is, $\alpha=(\alpha_1,\alpha_2,\ldots)\in {\mathbb N}^{\infty}$ having finitely many nonzero entries, where $\mathbb{N}=\{0,1,2,\ldots\}$.

\begin{Definition}[{A.~Lascoux~\cite{Las03}}]
The $\beta$-\emph{Lascoux polynomials} $\Omega_{\alpha}^{(\beta)}$ are again defined recursively.
For $\alpha\in {\sf Comp}$, define
\[\Omega^{(\beta)}_{\alpha}=
\begin{cases}
x_1^{\alpha_1} x_2^{\alpha_2} x_{3}^{\alpha_3}\cdots 
&\text{\   if $\alpha_1\geq \alpha_2\geq\alpha_3\geq\ldots$},\\
{\widetilde \pi}_i  (\Omega^{(\beta)}_{\alpha s_i} )
&\text{ \ if $\alpha_i<\alpha_{i+1}$.}
\end{cases}
\]
\end{Definition}

The nomenclature ``Lascoux polynomial'' first appears in C.~Monical's \cite{Monical:skyline}. We also refer to A.~Lascoux's \cite{Lascoux},
A.~N.~Kirillov's \cite{Kirillov}, C.~Monical-O.~Pechenik-D.~Searles' \cite{MPS}, O.~Pechenik-D.~Searles' \cite{PS19}, and the references therein for more about both families of polynomials.

\subsection{Increasing tableaux and $K$-jeu de taquin}\label{sec:Kjdt}

We need some notions from \cite{Thomas.Yong:K}.
\begin{Definition}
An \emph{increasing tableaux} of shape $\nu/\lambda$ is a filling of $\nu/\lambda$ using $\{1,2,\ldots,|\nu/\lambda|\}$ such that the
labels of $T$ strictly increase along rows and columns.
\end{Definition}
Let ${\sf INC}(\nu/\lambda)$ denote the set of all increasing tableaux of shape $\nu/\lambda$.

\begin{Definition}
A \emph{short ribbon} $R$ is a skew shape without a $2\times 2$ subshape, where
each row and column has at most two boxes, and each box is filled
with one of two symbols, but
adjacent boxes are filled differently.  Two boxes lie in the same \emph{component} of $R$ if there is a
path between them passing through boxes that are adjacent vertically or horizontally.
\end{Definition}

\begin{Definition}
Define
${\tt switch}(R)$ to be the same short ribbon as $R$ but where,
in each non-singleton component, each box is filled with the other symbol.
\end{Definition}

For example:
\[R=\tableau{&&&{\circ}\\&{\circ}&{\bullet}\\{\circ}&{\bullet}}
\mbox{ \ \ \ \ \ \ \ \ \ \ \ \ \ \
${\tt switch}(R)=
\tableau{&&&{\circ}\\&{\bullet}&{\circ}\\{\bullet}&{\circ}}$.}
\]
In what follows, we assume $\nu/\lambda$ is contained in an ambient rectangle $\Lambda$.

\begin{Definition}
An \emph{outer corner} of a skew shape $\nu/\lambda$ is a maximally northwest box of $\Lambda/\nu$.

\end{Definition}
Given $T\in {\sf INC}(\nu/\lambda)$, consider a \emph{set} of outer corners
$\{x_i\}$ filled with $\bullet$. Let $m$ be the maximum value label appearing in $T$.
Define ${\tt revKjdt}_{\{x_i\}}(T)$ as follows: let $R_m$ be the short ribbon consisting of $\bullet$ and $m$.
Apply ${\tt switch}(R_m)$. Now let $R_{m-1}$ be the short ribbon consisting of
$\bullet$ and $m-1$ and apply ${\tt switch}(R_{m-1})$.  Repeat until one applies
${\tt switch}(R_1)$, and then erase the $\bullet$ entries. For example, if $\nu/\lambda=(3,2,1)/(2,1)$ is contained in $\Lambda=(3,3,3,3)$ and we might have
\[T=\tableau{&&{2}\\&{2}&\bullet \\{1}&\bullet\\ \bullet}
\mapsto
\tableau{&&\bullet\\&{\bullet}&{2} \\{1}&2 \\ \bullet}
\mapsto
\tableau{&&\bullet\\&{\bullet}&{2} \\{\bullet }&2 \\ 1}
=\tableau{&&\\&&2\\&2\\1}={\tt revKjdt}_{\{x_i\}}(T).
\]
\begin{Definition}

\emph{A reverse $K$-rectification} of $T$ is any sequence of ${\tt revKjdt}$-slides
giving a reverse straight tableau ${\tt revKrect}(T)$.
\end{Definition}

Continuing the previous example, one can perform the following ${\tt revKjdt}$-slides:
\[\tableau{&&\\&&2\\&2&{\bullet}\\1} \ \ \
\tableau{&&\\&&\bullet\\&\bullet&{2}\\1} \ \ \
\tableau{&&\\&&\\&&{2}\\1&\bullet} \  \ \ 
\tableau{&&\\&&\\&&{2}\\\bullet&1} \  \ \ 
\tableau{&&\\&&\\&&{2}\\ &1&\bullet} \ \ \ 
\tableau{&&\\&&\\&&{\bullet}\\ &1&2}  \  \ \
\tableau{&&\\&&\\&&\\ &1&2}
\]
to conclude ${\tt revKrect}(T)=\tableau{1&2}$. Unlike textbook jeu de taquin, ${\tt revKrect}$ might 
depend on the choice
of ${\tt revKjdt}$-slides used (see \cite[Example~1.6]{Thomas.Yong:K}, \cite[Example~3.4]{Buch.Samuel}).

Given $P\in {\sf INC}(\nu)$ we define the \emph{left key} $K_{-}(P)$ to be a tableau, using 
${\tt revKjdt}$, as follows. By definition, the first columns of $P$ and $K_{-}(P)$ agree.  Assume that the first $\ell$ columns of $K_{-}(P)$ have been determined. Apply reverse rectification of the increasing tableau $P^{(\ell+1)}$ comprised of the first $\ell+1$ columns of $P$, inside the smallest rectangle $\Lambda^{(\ell+1)}$ that $P^{(\ell+1)}$ fits inside. \emph{For specificity, we define
${\tt revKrect}$ by using the leftmost outer corner 
for each intermediate ${\tt revKjdt}$-slide.}
Let $C^{(\ell+1)}$ be the leftmost column in the reverse rectification of $P^{(\ell+1)}$.
Then $C^{(\ell+1)}$ (after upward-justification) is the $(\ell+1)$-st column
of $K_{-}(P)$. Repeating this,  the end result  is $K_{-}(P)$.

\begin{Example}
\label{left-key-example}
The reader can check that if  
\[P=\tableau{1 & 2 & 3 & 5 & 7\\ 2 & 4 & 5 & 6\\ 4 & 6 } \text{ \ then  \ }
K_{-}(P)=\tableau{1 & 1 & 1 & 1 & 2 \\ 2 & 2 & 2 & 2\\ 4 & 4}
.\]
The first two columns are easily seen from the definition, as reverse rectification does nothing. To compute the third column
one works out this ${\tt revKrect}$:
\[\tableau{1 & 2 &3\\ 2 & 4 & 5 \\ 4 &6 &\bullet}\to 
\tableau{1 & 2 &3\\ 2 & 4 & 5 \\ 4 &\bullet  & 6 }\to
\tableau{1 & 2 &3\\ 2 & \bullet & 5 \\ \bullet &4 & 6}\to
\tableau{1 & \bullet &3\\ \bullet & 2 & 5 \\ 2 &4 & 6}\to
\tableau{\bullet & 1 &3\\ 1 & 2 & 5 \\ 2 &4 & 6}.
\]
The first column $C^{(3)}=\tableau{1 \\ 2}$ of the rightmost tableau is the third column of
$K_{-}(P)$. \qed
\end{Example}

Notice in our example, $K_{-}(P)$ has the same (straight) shape as $P$. 
Also, $K_{-}(P)$ is a \emph{key}: the set of labels in column $i$ are contained in the
set of labels in column $i-1$ for $i\geq 2$.  These properties always hold, and are proven 
in Section~\ref{sec:aproof}.  Let ${\sf content}(K_{-}(P))$ be the usual content of a semistandard
tableau; here ${\sf content}(K_{-}(P))=(4,5,0,2)$.

\subsection{Reduced and Hecke words}
Let $\ell(w)$ be the \emph{Coxeter length} of $w\in S_n$, that is, 
\[\ell(w)=\#\{1\leq i<j\leq n: w(i)>w(j)\}.\]

\begin{Definition}
A sequence $(i_1,i_2,\ldots,i_{\ell(w)})$ is a \emph{reduced word} for $w$ if $s_{i_1}s_{i_2}\cdots s_{i_{\ell(w)}}=w$.
\end{Definition}
Let ${\sf Red}(w)$ denote the set of reduced words for $w$.

\begin{Definition}
The \emph{Hecke monoid} ${\mathcal H}_{n}$ is generated by $u_1,u_2,\ldots,u_{n-1}$, subject to:
\begin{align*}
u_i^2 & \equiv   u_i\\
u_i u_j & \equiv u_j  u_i \text{\ \ \ \  if $|i-j|>1$}\\
u_i  u_{i+1}  u_i & \equiv u_{i+1} u_i  u_{i+1}
\end{align*} 
\end{Definition}
\begin{Definition}
A sequence $(i_1,i_2,\ldots,i_N) \in {\mathbb N}^N$ is a \emph{Hecke word} for $w\in S_n$ if
\[
u_{i_1} u_{i_2}\cdots u_{i_N}\equiv  u_{a_1} u_{a_2}\cdots u_{a_{\ell(w)}},
\text{\  for some $(a_1,\ldots, a_{\ell(w)})\in {\sf Red}(w)$.}
\]
\end{Definition}

\begin{Definition}
For any tableau $P$, we will read off a \emph{word} denoted ${\sf word}(P)$, concatenating its rightmost column read top-to-bottom, then its next-to-rightmost column, etc.  For example, the tableau $P$ in Example~\ref{left-key-example} has ${\sf word}(P)=(7,5,6,3,5,2,4,6,1,2,4)$.
\end{Definition}

\subsection{The ``Grothendieck to Lascoux'' conjecture} 
This is the open problem of this report:

\begin{Conjecture}
\label{conj:main}
\[{\mathfrak G}^{(\beta)}_w=\sum_{P} \beta^{\# {\sf boxes}(P)-\ell(w)}\Omega^{(\beta)}_{{\sf content}(K_{-}(P))}\]
where $P$ is any straight-shape increasing tableau such that
${\sf word}(P)$ is a Hecke word for $w$.
\end{Conjecture}

\begin{Example}
If $w=31524$ the increasing tableaux and the left keys are respectively 
\begin{equation}
\label{eqn:Nov2abc}
P = \tableau{1 & 2 &4 \\ 3}, \ \ \ \tableau{1 & 2\\ 3 &4}, \ \ \ \tableau{1 & 2 & 4 \\ 3 &4}; \ \ 
K_{-}(P)=\tableau{1 & 1 &1 \\ 3}, \ \ \ \tableau{1 & 1\\ 3 & 3}, \ \ \ \tableau{1 & 1 & 1 \\ 3 &3}.
\end{equation} 
For instance, if $P$ is the rightmost increasing tableaux, then ${\sf word}(P)=(4,2,4,1,3)$. Now
\[u_4 u_2 u_4 u_1 u_3\equiv u_2 u_4^2 u_1 u_3\equiv u_2 u_4 u_1 u_3,\] 
and $s_2 s_4 s_1 s_3 =31524$.  Conjecture~\ref{conj:main} predicts 
${\mathfrak G}^{(\beta)}_{31524}=\Omega_{301}^{(\beta)}+\Omega_{202}^{(\beta)} + \beta \, \Omega_{302}^{(\beta)}$.\qed
\end{Example}

Conjecture~\ref{conj:main} generalizes one formula \cite{Reiner.Shimo} and refines another
\cite{BKSTY}. It has been exhaustively checked (with computer assistance) for $n\leq 7$ (and spot-checked for $n=8,9$). 
It says the  $\Groth_w^{(\beta)}$ to
$\Omega^{(\beta)}_{\alpha}$ expansion is positive; this  is also open.
%
%
%
 
The rest of this report surveys known and related 
results that motivate Conjecture~\ref{conj:main}. The proof that the shapes of $P$ and 
 $K_{-}(P)$ agree, and that $K_{-}(P)$ is a key, is in Section~\ref{sec:aproof}.

\section{History of the problem}\label{sec:2}

During the preparation of \cite{BKSTY}, M.~Shimozono privately conjectured to the second author that 
${\mathfrak G}_w^{(\beta)}$ expands
positively in the $\Omega_{\alpha}^{(\beta)}$'s; he also suggested ideas towards a rule. 
Conjecture~\ref{conj:main} was formulated in September 2011 during a visit of the first author to UIUC. There are two limiting cases of Conjecture~\ref{conj:main}, as explained now.

\subsection{The $\beta=0$ specialization and stable-limit symmetry}
\begin{Definition}
The \emph{key polynomial} (or type $A$ \emph{Demazure character}) is $\kappa_{\alpha}:=\Omega_{\alpha}^{(0)}$. 
\end{Definition}
References about key polynomials include \cite{LasSch90, Reiner.Shimo, Lascoux}. A tableau formula for $\kappa_{\alpha}$ is in \cite{LasSch90} (see also \cite{Reiner.Shimo}).  From the definitions, 
\begin{equation}
\label{eqn:deformationNov3}
\Omega_{\alpha}^{(\beta)}=\kappa_{\alpha}+\sum_{k>0}\beta^k p_k,
\end{equation} 
where $p_k$ is a homogeneous polynomial in $x_1,x_2,\ldots$
of degree $|\alpha|+k$, where $|\alpha|=\sum_{i\geq 1}\alpha_i$. 

\begin{Definition}\label{def:Schubert}
The \emph{Schubert polynomial}\footnote{It represents the class of a Schubert
variety $X_w$ in the flag variety $GL_n/B$ under Borel's isomorphism (see \cite{Fulton}). The \emph{Grothendieck polynomial} 
${\mathfrak G}_w:={\mathfrak G}_w^{(-1)}$ similarly represents the Schubert structure sheaf ${\mathcal O}_{X_w}$
in the Grothendieck ring $K^0(GL_n/B)$ of algebraic vector bundles on $GL_n/B$. This explains the ``combinatorial
$K$-theory'' nomenclature \cite{Buch:CKT}.}  is ${\mathfrak S}_w:={\mathfrak G}^{(0)}_w$. 
\end{Definition}
\vspace{-.1in}

A combinatorial rule for ${\mathfrak G}_w^{(\beta)}$ as a sum of ${\mathfrak S}_v$'s is given by C.~Lenart's \cite{Lenart99}.

All the aforementioned polynomial families (Lascoux, Key, Grothendieck, Schubert) are ${\mathbb Z}[\beta]$-linear
bases of ${\mathbb Z}[\beta][x_1,x_2,\ldots]$.  There are symmetric versions of these polynomials.    
\begin{Definition}
The $\beta$-\emph{stable Grothendieck polynomial} is
\[G^{(\beta)}_w(x_1,x_2,\ldots)=\lim_{n\to\infty}{\mathfrak G}^{(\beta)}_{1^n\times w}(X),\]
where $(1^n\times w)(i)=i$ if $1\leq i\leq n$ and $(1^n\times w)(i)=w(i)+n$ if $i>n$.
\end{Definition}

\begin{Definition}\label{def:Stanley}
The \emph{stable Schubert polynomial} is $F_w:=G_w^{(0)}$. This is also known as the \emph{Stanley symmetric polynomial}.
\end{Definition}

\begin{Definition}
A permutation $w$ is \emph{Grassmannian at position $k$} if $w(i)< w(i+1)$
for  $i\neq k$. To such $w$, define a partition $\lambda=\lambda(w)$ by $\lambda_i=w(k-i+1)-(k-i+1)$ for $1\leq i\leq k$. 
\end{Definition}

\begin{Definition}
A \emph{set-valued semistandard Young tableaux} $T$ of shape $\lambda$ is a filling of the boxes of $\lambda$ with nonempty sets such that if one chooses a singleton from each set, 
the result is a semistandard Young tableaux (row weakly increasing and column strict). 
\end{Definition}

\begin{Theorem}[\cite{Buch}]
Let $w$ be a Grassmannian permutation of shape $\lambda$.  Then
\[
G^{(\beta)}_{\lambda}:=G^{(\beta)}_w=\sum_{T}\beta^{\#{\sf labels}(T)-|\lambda|}x^T,
\]
where the sum is over set-valued semistandard Young tableaux $T$ of shape $\lambda$,
and $x^T:=\prod_{i\geq 1} x_i^{\#i\in T}$.
\end{Theorem}

\begin{Definition}
The \emph{Schur function} is $s_{\lambda}:=G^{(0)}_{\lambda}$. 
\end{Definition}

Summarizing, one has a commutative diagram
\begin{equation}
\label{KDemazure-square}
\begin{array}{ccc}
\KDemazure_\alpha^{(\beta)} & \longrightarrow & \Demazure_\alpha \\
\downarrow & & \downarrow \\
G_{\lambda(\alpha)}^{(\beta)} & \longrightarrow & s_{\lambda(\alpha)} \\
\end{array}
\end{equation}
with horizontal arrows indicating $\beta=0$ specialization, and vertical arrows called \emph{stabilization}: for
the right vertical arrow, let $\lambda(\alpha)$ be the sorting of $\alpha$, then if $N>n$, 
\[s_{\lambda(\alpha)}(x_1,\ldots,x_n)=\kappa_{(0^N,\alpha)}(x_1,\ldots,x_n,0,0,\ldots \ ).\]

\subsection{Monomial expansion formulas}

\begin{Theorem}[\cite{Fomin.Kirillov}]\label{grothform}
\[{\mathfrak G}^{(\beta)}_{w}(X)=\sum_{({\bf a},{\bf i})} \beta^{N-\ell(w)}x^{\bf i}\]
where the sum is over all pairs of sequences $({\bf a},{\bf i})$ (called compatible sequences) such that
\begin{itemize}
\item[(a)] ${\bf a}=(a_1,a_2,\cdots,a_N)$ is a Hecke word for $w$;
\item[(b)] ${\bf i}=(i_1,i_2,\cdots,i_N)$ has $1\leq i_1\leq i_2\leq\ldots\leq i_N$
\item[(c)] $i_j\leq a_j$; and
\item[(d)] $a_j\leq a_{j+1}\implies i_j<i_{j+1}$.
\end{itemize}
\end{Theorem}

\begin{Theorem}[\cite{Fomin.Kirillov}]\label{stablegrothform}
\[G_w^{(\beta)}=\sum_{({\bf a},{\bf i})} \beta^{N-\ell(w)}x^{\bf i}\]
where $({\bf a},{\bf i})$ satisfies (a), (b) and (d) above.
\end{Theorem}

Therefore, one has a commutative diagram,
with the same arrows
\begin{equation}
\label{Grothendieck-square}
\begin{array}{ccc}
\Groth_w^{(\beta)}\overset{{\rm Thm~\ref{grothform}}}{=}\sum_{({\bf a},{\bf i})} \beta^{N-\ell(w)} x^{\bf i}& \longrightarrow & \Schub_w=\sum_{({\bf a},{\bf i})} x^{\bf i}\\
\downarrow & & \downarrow \\
G_w^{(\beta)}\overset{{\rm Thm~\ref{stablegrothform}}}{=}\sum_{({\bf a},{\bf i})} \beta^{N-\ell(w)} x^{\bf i} & \longrightarrow & F_w=\sum_{({\bf a},{\bf i})} x^{\bf i}. \\
\end{array}
\end{equation}

In view of Definitions~\ref{def:Schubert} and~\ref{def:Stanley}, the expressions in the right column sum over
$({\bf a},{\bf i})$ with ${\bf a} \in {\sf Red}(w)$, versus sums over Hecke words ${\bf a}$ for $w$ in their left column counterparts.

The diagram \eqref{Grothendieck-square} specializes when $w=w(\lambda)$ is Grassmannian, 
giving
\begin{equation}
\label{Grassmannian-square}
\begin{array}{ccc}
G^{(\beta)}_\lambda(x_1,\ldots,x_n) & \longrightarrow & s_\lambda(x_1,\ldots,x_n)\\
\downarrow & & \downarrow \\
G^{(\beta)}_\lambda & \longrightarrow & s_\lambda. \\
\end{array}
\end{equation}

\subsection{Prior expansion formulas}
Conjecture~\ref{conj:main} generalizes a relationship between the Schubert and key polynomials.

\begin{Theorem}[\cite{LS89, Reiner.Shimo}]\label{thm:Schubtokey}
\[{\mathfrak S}_w=\sum_{P} \kappa_{{\sf content}(K_{-}(P))}\]
where the sum is over all increasing tableaux $P$ such that ${\sf word}(P)\in {\sf Red}(w)$.
\end{Theorem}
To be precise, in the formulation given in \cite[Theorem~4]{LS89}, the description of the ``left nil key''
$K_{-}(P)$ differs from our definition. We are asserting (proof omitted) that in the case of the $P$ in Theorem~\ref{thm:Schubtokey},
the two definitions agree. This is because ${\tt revKjdt}$ can be used to compute the insertion tableau of \emph{Hecke
insertion} \cite{BKSTY} which specializes to \emph{Edelman-Greene insertion} \cite{EG}; see \cite{TY:Advapp}.

Theorem~\ref{thm:Schubtokey} is the non-symmetric version of the following result:
\begin{Theorem}[\cite{fomin.greene:noncommutative}]\label{thm:FGEG}
Let $a_{w,\lambda}=\#\{P\in {\sf INC}(\lambda): {\sf word}(P)\in {\sf Red}(w)\}$. Then
\[F_w=\sum_{\lambda} a_{w,\lambda}s_{\lambda}.\]
\end{Theorem}

The next result generalizes Theorem~\ref{thm:FGEG}.  
Conjecture~\ref{conj:main} is
its non-symmetric version:
\begin{Theorem}[\cite{BKSTY}]\label{thm:BKSTY}
Let $b_{w,\lambda}=\#\{P\in {\sf INC}(\lambda)\!:\! {\sf word}(P)\text{ is a Hecke word for }w\}$. Then
\[G^{(\beta)}_{w}=\sum_{\lambda} \beta^{|\lambda|-\ell(w)} \ b_{w,\lambda} G_{\lambda}^{(\beta)}.\]
\end{Theorem}

\begin{Example}\label{exa:31524}
$G_{31524}^{(\beta)}=G_{31}^{(\beta)}+G_{22}^{(\beta)}+\beta G_{32}^{(\beta)}$. This is witnessed by the $P$ tableaux of
(\ref{eqn:Nov2abc}). Notice that Conjecture~\ref{conj:main} subdivides the witnessing tableaux $P$ for $b_{w,\lambda}$
according to ${\sf content}(K_{-}(P))$. It is in this sense that Conjecture~\ref{conj:main} is a refinement of Theorem~\ref{thm:BKSTY}.\qed
\end{Example}

In conclusion, Conjecture~\ref{conj:main}
captures some known facts, as expressed in
this diagram
\begin{equation}
\label{tableau-expansion-square}
\begin{array}{ccc}
\Groth_w^{(\beta)} \overset{\mathrm{Conj~\ref{conj:main}}}{=}\sum_P \beta^{|{\sf shape}(P)|-\ell(w)} 
\KDemazure_{{\sf content}(K_-(P))}^{(\beta)} &
   \longrightarrow & \Schub_w \overset{\mathrm{Thm~\ref{thm:Schubtokey}}}{=} \sum_P \Demazure_{{\sf content}(K_-(P))} \\
\downarrow & & \downarrow \\
G_w^{(\beta)} \overset{\mathrm{Thm~\ref{thm:BKSTY}}}{=}\sum_P  \beta^{|{\sf shape}(P)|-\ell(w)} G_{{\sf shape}(P)}^{(\beta)} &
   \longrightarrow &F_w \overset{\mathrm{Thm~\ref{thm:FGEG}}}{=} \sum_P s_{{\sf shape}(P)}.\\
\end{array}
\end{equation}
Here ${\sf shape}(P)$ is the partition of $P$. The horizontal, vertical maps
are as before.

\section{Further discussion}

\subsection{Formulas for Lascoux polynomials}
\vspace{-.1in}
Combinatorial rules for the Lascoux polynomials are in V.~Buciumas--T.~Scrimshaw--K.~Weber
\cite{Scrimshaw}. Another rule, generalizing the Kohnert moves of \cite{Kohnert}, 
was conjectured in \cite{Ross.Yong}.\footnote{The conjecture is accidentally misstated there.
See the corrected version \url{https://faculty.math.illinois.edu/~ayong/polynomials.Seminaire.revision.2017.pdf} which is
consistent with the 2011 report by C.~Ross \url{https://faculty.math.illinois.edu/~ayong/student_projects/Ross.pdf}.} 
Also, see the skyline conjectural rule of C.~Monical \cite{Monical:skyline}. Finally,
\cite[Theorem~5]{Reiner.Shimo} gives an alternative formula for $\kappa_{\alpha}$ in terms of compatible sequences; we do
not know a generalization of this formula to $\Omega_{\alpha}^{(\beta)}$. 

\subsection{Warning about stable-limits} The results of Section~\ref{sec:2} suggest  
combinatorial properties for stable-limit polynomials will hold for their non-symmetric versions. This is not always true. S.~Fomin-C.~Greene proved the following result (\emph{cf.} C.~Lenart's \cite{Lenart00}):
\begin{Theorem}[\cite{fomin.greene:noncommutative}] \label{thm:FGthing}
\begin{equation}
\label{eqn:FGthingexp}
G_w^{(\beta)}=\sum_{\lambda} \beta^{|\lambda|-\ell(w)}d_{w,\lambda}s_{\lambda}
\end{equation}
where $d_{w,\lambda}$ counts
tableaux $P$ of shape $\lambda$ that are row strict and column weakly increasing,
such that ${\sf word}(P)$  is a Hecke word for $w$.
\end{Theorem}
Thus, using the Grassmannian permutation $w=s_4 s_1 s_2 s_3 = 23514$, 
\[G_{2,1,1}^{(\beta)}=s_{2,1,1}+\beta(3s_{2,1,1,1}+s_{2,2,1})+\cdots.\] 
If one expands $\Omega^{(\beta)}_{\alpha}$ in the keys (the 
non-symmetric analogue (\ref{eqn:FGthingexp}); see (\ref{KDemazure-square})),  
positivity fails:
\[\Omega^{(\beta)}_{1,0,2,1}=\kappa_{1,0,2,1}+\beta (2\kappa_{1,1,2,1}+\kappa_{2,0,2,1}+\kappa_{1,2,2}
-\kappa_{2,1,2})+\cdots.\]

\subsection{Proof that $K_{-}(P)$ is a key of the same shape as $P$}\label{sec:aproof}
We first show ${\sf shape}(K_{-}(P))={\sf shape}(P)$. 
In the notation of $K_{-}(P)$'s description,  it suffices to argue
that the length $b$ of $C^{(\ell+1)}$ equals the length $t$ of the $\ell+1$ (\emph{i.e.}, rightmost) column of $P^{(\ell+1)}$.

To see this, consider the general situation of an increasing tableau $T$ contained in an $r\times s$ dimension rectangle $\Lambda$, and a second
increasing tableau in $\Lambda$ that is complementary to $T$. Recall the notion of ${\tt Kinfusion}$ defined in \cite[Section~3]{Thomas.Yong:K}.
In fact, 
\begin{equation}
\label{eqn:Feb9abc}
{\tt Kinfusion}(T,U)=(A,B) \text{\ where $A={\tt Krect}(U)$ and $B={\tt revKrect}(A)$.}
\end{equation}
The rectification $A$ uses the
inner corners defined by $T$. Similarly the reverse rectification $B$ uses the outer corners defined by $U$.

Given $V\in {\sf INC}(\nu/\lambda)$ define ${\sf LDS}(V)$ to be the length of the longest decreasing subsequence of the left to right, bottom to top, row reading word of $V$. This is true:
\begin{Theorem}[{\cite[Theorem~6.1]{Thomas.Yong:K}, \emph{cf.} \cite[Corollary~6.8]{Buch.Samuel}}]\label{thm:LDS}
${\sf LDS}$ is invariant under $K$-theoretic
(reverse) jeu de taquin slides. 
\end{Theorem}

Now, suppose $t$ and $u$ are the lengths of the $s$-th (possibly empty) column of $T$ and $U$, respectively. 
Similarly, let $a$ and $b$ be the length of the first (leftmost) columns of $A$ and $B$, respectively. Thus $t+u=a+b=r$.
By Theorem~\ref{thm:LDS}, $u={\sf LDS}(U)={\sf LDS}(A)=a$. Thus $b=t$. The result follows from (\ref{eqn:Feb9abc})
and setting $T=P^{(\ell+1)}$ and $U$ being any complementary increasing tableau to $T$ inside the smallest rectangle $\Lambda^{(\ell+1)}$ that 
$P^{(\ell+1)}$ sits inside.

To see that $K_{-}(P)$ is a key we use an argument of G.~Orelowitz: Since we choose the leftmost
outer corner at each slide of ${\tt revKrect}$, when computing $C^{(\ell+1)}$ we
begin by computing ${\tt revKrect}(P^{(\ell)})$ as a partial ${\tt revKrect}$ of $P^{(\ell+1)}$. At this point, $C^{(\ell)}$ is the leftmost
column of this partial {\tt revKrect}. Thus, when completing the ${\tt revKrect}(P^{(\ell+1)})$, by the
 definition of ${\tt revKjdt}$, the entries of $C^{(\ell+1)}$ are contained in those of $C^{(\ell)}$, as desired.
\qed

The above shape argument does not depend on the specific choice of ${\tt revKrect}$ used at each stage of the definition of
$K_{-}(P)$. We suspect this choice does not affect $K_{-}(P)$ being a key, however, the choice we use (suggested by G.~Orelowitz) 
makes the proof easy.

\section*{Acknowledgements}
AY thanks Mark Shimozono for initiating his interest in the expansion problem of this report.
We are grateful to Shiliang Gao, Cara Monical, and Gidon Orelowitz who each gave very helpful analysis. We thank Oliver Pechenik, Colleen Robichaux, Colleen Ross, Travis Scrimshaw, and Dominic Searles for discussions. AY is supported by a Simons Collaboration Grant, NSF RTG grant DMS-1937241, and the UIUC
Center for Advanced Study. VR is supported by NSF grant DMS-1601961. John Stembridge's Maple package {\sf Coxeter} was used in our experiments. We used computing at Brown's 
Center for Computation and Visualization during AY's virtual residence at ICERM in Spring 2021.


\end{document}